# A note on Turán number $T(k+1,k,n)$

Li An-Ping

Beijing 100085, P.R. China
apli0001@sina.com

Abstract: Turán number is one of primary topics in the combinatorics of finite sets, in this paper, we will present a new upper bound for Turán number $T(k+1,k,n)$.

## 1. Introduction

Suppose that $S$ is a set of $n$ elements, for a non-negative integer $r$, as usual, $\binom{S}{r}$ stand for the set of all the $r-subsets$ of $S$. Suppose that $k$ is a non-negative integer, $k \leq r$, $\mathcal{F} \subseteq \binom{S}{k}$ is a family of $k-subsets$ of $S$, $\mathcal{F}$ is called a $(r,k)-$ Turán family if for any $X \in \binom{S}{r}$, there is a member $A \in \mathcal{F}$ such that $A \subseteq X$. Denoted by $T(r,k,n)$ the minimal size of $(r,k)-$ Turán families, the number is also called Turán number, which first proposed in P. Turán's papers [5], [6].

For $T(r,k,n)$, there are the estimations

$$T(r,k,n) \geq \binom{n}{k} \Big/ \binom{r}{k}, \qquad (1.1)$$

and

$$T(r,k,n) \leq \sum_{i=1}^{d} \binom{|S_i|}{k} \qquad (1.2)$$

where $d = \left[\frac{r-1}{k-1}\right]$, $|S_i| = \left[\frac{n}{d}\right]$, or $|S_i| = \left\lceil\frac{n}{d}\right\rceil$, so for $k > 3$

$$T(r,k,n) \leq \frac{1}{d^{k-1}} \binom{n}{k}. \qquad (1.3)$$

The details for the results above refer to see papers [2], [4]~[6].

The result (1.3) will be trivial in the case $d=1$, i.e. the case $r < 2k-1$, for instance, $r = k+1$.

For $T(k+1,k,n)$, K.H. Kim and F.W. Roush [3] proved that

$$T(k+1,k,n) \leq \left(\frac{1}{[k/2\log k]} + \frac{1}{2k\log k}\right)\binom{n}{k}, \qquad (1.4)$$

V.Rödl and P. Frankl [1] improved that

$$T(k+1,k,n) \leq \frac{1}{[k/\log k]} \frac{\log k}{\log k - 1} \binom{n}{k}. \qquad (1.5)$$

In this paper, our main result is that.

**Theorem 1.**

$$T(k+1,k,n) \le \frac{\log k + c_1(k) - c_0(k)}{k-1} C_n^k, \qquad (1.6)$$

where $c_0(k) = \int_1^{\log k} (1-e^{-x})^k dx$, $c_1(k) = \sum_{i>0}(-1)^{i-1} C_k^i k^{-i}/i$, $c_1(k) \le 1$, $c_0(k) \ge 0$.

## 2. Proof of Theorem 1

As usual, $\mathbb{Z}_0$ stand for the set of all non-negative integers, suppose that $s$, $k$ and $h$ are three positive integers, denoted by

$$\Omega_s = \left\{ v \mid v = (x_1, x_2, \ldots, x_k) \in \mathbb{Z}_0^k, \sum_{1 \le i \le k} x_i = s \right\},$$

$$V(s,k,h) = \{v \mid v \in \Omega_s, x_i \le h, 1 \le i \le k\}, \quad \overline{V}(s,k,h) = \left\{v \mid v \in \Omega_s, \max_{1 \le i \le k}\{x_i\} = h\right\}. \qquad (2.1)$$

Define

$$\mu(s,k,h) = \sum_{(x_1,\ldots,x_k) \in V(s,k,h)} (x_k + 1), \quad \lambda(s,k,h) = \sum_{(x_1,\ldots,x_k) \in \overline{V}(s,k,h)} (x_k + 1). \qquad (2.2)$$

**Lemma 1.**

$$\mu(s,k,h) = \sum_i (-1)^i \left( C_k^i C_{s+k-i(h+1)}^k + (h+1) C_{k-1}^{i-1} C_{s+k-i(h+1)-1}^{k-1} \right). \qquad (2.3)$$

*Proof.* Suppose that $g(x) = \sum_s \mu(s,k,h) \cdot x^s$ is the generating function of $\mu(s,k,h)$ for the parameter $s$. Denoted by $h_1 = h+1$, let $p(x) = (1-x^{h_1})/(1-x)$, $p_1(x) = (1-x^{h_1+1})/(1-x)$, then it has

$$g(x) = p^{k-1}(x) \times p_1(x)'$$
$$= \left(1-x^{h_1}\right)^{k-1} \times \left(1-(h_1+1)x^h + h_1 x^{h+1}\right) \times \sum_j C_{k+j}^k x^j$$
$$= \sum_i (-1)^i (C_k^i - h_1 x C_{k-1}^{i-1} + h_1 C_{k-1}^{i-1}) x^{i(h+1)} \times \sum_j C_{k+j}^k x^j$$
$$= \sum_s x^s \sum_i (-1)^i \left( C_k^i C_{s+k-ih_1}^k + h_1 C_{k-1}^{i-1} C_{s+k-ih_1-1}^{k-1} \right)$$

i.e.

$$\mu(s,k,h) = \sum_i (-1)^i \left( C_k^i C_{s+k-ih_1}^k + h_1 C_{k-1}^{i-1} C_{s+k-ih_1-1}^{k-1} \right). \qquad \square$$

Let $s = n-k$, and assumed that $n \gg k$.

**Lemma 2.**

$$\sum_{0 \leq j < k} \log\left(1 - \frac{v}{n-j}\right) = v\log\left(1 - \frac{k}{n}\right) - \frac{1}{2}kv^2 n^{-2} + o\left(kv^2 n^{-2}\right). \quad (2.4)$$

*Proof.* From the calculus we know

$$\int_0^k \log\left(1 - \frac{v}{n-x}\right)dx \leq \sum_{0 \leq j < k} \log\left(1 - \frac{v}{n-j}\right) \leq \int_0^k \log\left(1 - \frac{v}{n-x}\right)dx - \frac{1}{2}\left(\log\left(1 - \frac{v}{n-k}\right) - \log\left(1 - \frac{v}{n}\right)\right)$$

and

$$\int_0^k \log\left(1 - \frac{v}{n-x}\right)dx = -\int_0^k \sum_{i>0} \frac{1}{i}\left(\frac{v}{n-x}\right)^i dx$$

$$= v\log\left(1 - \frac{k}{n}\right) - v \cdot \sum_{i>0}\left(\left(\frac{v}{n-k}\right)^i - \left(\frac{v}{n}\right)^i\right) \cdot \frac{1}{i(i+1)} = v\log\left(1 - \frac{k}{n}\right) - \sum_{i \geq 2}\left(\frac{v}{n-\theta}\right)^i \cdot \frac{k}{i}$$

$$= v\log\left(1 - \frac{k}{n}\right) - \frac{1}{2}kv^2 n^{-2} + o(kv^2 n^{-2}). \qquad \square$$

**Lemma 3.**

$$\mu(s,k,h-1) \geq C_n^k \times \left(\left(1 - (1-(k/n))^h\right)^d - \frac{hk}{n}\left(1 - (1-(k/n))^h\right)^{k-1} \cdot \left(1-(k/n)\right)^h\right)$$
$$-(1/2 + o(1)) \cdot C_n^k \cdot k \cdot \sum_{i \geq 1}(-1)^i C_k^i (1-(k/n))^{ih}((hi/n)^2 + (hi/n)^3) \quad (2.5)$$

*Proof.* By (2.3) and (2.4), it has

$$\mu(s,k,h-1) = \sum_i (-1)^i (C_k^i \cdot C_{n-ih}^k + hC_{k-1}^{i-1} \cdot C_{n-ih-1}^{k-1})$$

$$= C_n^k \times \sum_i (-1)^i \left(C_k^i \cdot \prod_{0 \leq j < k}\left(1 - \frac{i \cdot h}{n-j}\right) + \frac{hk}{n} C_{k-1}^{i-1} \cdot \prod_{0 < j < k}\left(1 - \frac{i \cdot h}{n-j}\right)\right)$$

$$\geq C_n^k \times \left(\left(1 - (1-(k/n))^h\right)^k - \frac{hk}{n}\left(1 - (1-((k-1)/(n-1)))^h\right)^{k-1} \cdot \left(1-((k-1)/(n-1))\right)^h\right)$$

$$-(1/2 + o(1)) \cdot C_n^k \cdot k \cdot \sum_{i \geq 1}(-1)^i C_k^i (1-(k/n))^{ih}((hi/n)^2 + (hi/n)^3)$$

$$\geq C_n^k \times \left(\left(1 - (1-(k/n))^h\right)^k - \frac{hk}{n}\left(1 - (1-(k/n))^h\right)^{k-1} \cdot \left(1-(k/n)\right)^h\right)$$

$$-(1/2 + o(1)) \cdot C_n^k \cdot k \cdot \sum_{i \geq 1}(-1)^i C_k^i (1-(k/n))^{ih}((hi/n)^2 + (hi/n)^3).$$

$\square$

**Lemma 4.** Let $\tau(t) = \int_1^t 1-(1-e^{-x})^t dx$, then

$$\tau(t) \leq \log(t) - 1 - c_0(t) + c_1(t). \qquad (2.6)$$

where $c_0(t) = \int_1^{\log t} (1-e^{-x})^t dx$, $c_1(t) = \sum_{i>0} (-1)^{i-1} C_t^i t^{-i} / i$, $c_1(t) \leq 1$, for $t \geq 1$.

And let $\tau_1(t) = \int_1^t (1-e^{-x})^{t-1} xe^{-x} dx$, then

$$\tau_1(t) \leq \frac{1+\tau(t)}{t}. \qquad (2.7)$$

*Proof.*

$$\int_1^t 1-(1-e^{-x})^t dx = \int_1^{\log(t)} 1-(1-e^{-x})^t dx + \int_{\log(t)}^t 1-(1-e^{-x})^t dx$$

$$\leq \log(t) - 1 - c_0(t) + \sum_{i>0} (-1)^{i-1} C_t^i t^{-i} / i.$$

And

$$\tau(t) = (t-1) - (1-e^{-x})^t x \big|_1^t + \int_1^t t(1-e^{-x})^{t-1} xe^{-x} dx$$

$$= (t-1) - (1-e^{-t})^t t + (1-e^{-1})^t + t \cdot \tau_1(t)$$

$$\geq t \cdot \tau_1(t) - 1. \qquad \square$$

**Lemma 5.**

$$\sum_{\lceil s/k \rceil \leq h < n-k} \mu(s,k,h) \geq \left( n - k - \frac{n}{k-1}(1+\tau(k)) \right) \cdot C_n^k. \qquad (2.8)$$

*Proof.* With (2.5), (2.6) and (2.7), it has

$$\sum_{\lceil s/k \rceil \leq h < n-k} \mu(s,k,h)$$

$$\geq C_n^k \cdot \sum_{\lceil n/k \rceil \leq h \leq n-k} \left( \left(1-(1-(k/n))^h\right)^d - \frac{hk}{n}\left(1-(1-(k/n))^h\right)^{k-1} \cdot (1-(k/n))^h \right)$$

$$- (1/2 + o(1)) \cdot C_n^k \cdot k \sum_{\lceil n/k \rceil \leq h \leq n-k} \sum_{i \geq 1} (-1)^i C_k^i (1-(k/n))^{ih} ((hi/n)^2 + (hi/n)^3)$$

$$\geq C_n^k \cdot \left( \sum_{\lceil n/k \rceil \leq h \leq n-k} (1-e^{-k \cdot h/n})^k - \sum_{\lceil n/d \rceil \leq h \leq n-k} \frac{hk}{n}(1-e^{-k \cdot h/n}))^{k-1} \cdot e^{-k \cdot h/n} \right)$$

$$-(1/2+o(1))\cdot C_n^k \cdot k \sum_{i\geq 1}(-1)^i C_k^i \sum_{\lceil n/k \rceil \leq h \leq n-k} (1-(k/n))^{ih}((hi/n)^2+(hi/n)^3)$$

$$\geq C_n^k \cdot \left( \int_{\lceil n/k \rceil -1}^{n-k} (1-e^{-kx/n})^k dx - \int_{\lceil n/k \rceil -1}^{n-k} (1-e^{-kx/n})^{k-1} \frac{k}{n} xe^{-kx/n} dx \right)$$

$$-(1+o(1))\cdot C_n^k \cdot \left(\frac{n}{k^2}+\frac{3n}{k^3}\right) \cdot \sum_{i\geq 1}(-1)^i C_k^i e^{-i}$$

$$\geq C_n^k \times \frac{n}{k} \times \left( \int_1^{(n-k)k/n} (1-e^{-x})^k dx - \int_1^{(n-k)k/n} (1-e^{-x})^{k-1} xe^{-x} dx \right)$$

$$+(1+o(1))\cdot C_n^k \cdot \left(\frac{n}{k^2}+\frac{3n}{k^3}\right) \cdot (1-(1-e^{-1})^k)$$

$$\geq C_n^k \times \left( n-k - \frac{n}{k} \times \frac{k+1}{k}(1+\tau(k)) \right) \geq C_n^k \times \left( n-k - \frac{n}{k-1}(1+\tau(k)) \right).$$

□

**Lemma 6.**

$$\sum_{\lceil s/k \rceil \leq h \leq n-k} h \cdot \lambda(s,k,h) \leq \left( \frac{n}{k-1}(1+\tau(k)) \right) \cdot C_n^k. \tag{2.9}$$

*Proof.* It is clear that

$$\lambda(s,d,h) = \mu(s,d,h) - \mu(s,d,h-1).$$

So,

$$\sum_{\lceil s/k \rceil \leq h \leq n-k} h \cdot \lambda(s,k,h) = \sum_{\lceil s/k \rceil \leq h \leq n-k} h \cdot (\mu(s,k,h) - \mu(s,k,h-1))$$

$$\leq (n-k)\cdot C_n^k - \sum_{\lceil s/k \rceil \leq h < (n-k)} \mu(s,k,h)$$

$$\leq \left( \frac{n}{k-1}(1+\tau(k)) \right) \cdot C_n^k. \qquad \square$$

*The Proof of Theorem 1.*

For two non-negative integers $a$ and $b$, $a<b$, denoted by $[a,b]=\{a,a+1,\ldots,b\}$, $(a,b)=\{a+1,a+2,\ldots,b-1\}$, so, in the case $b=a+1$, $(a,b)=\varnothing$.

Assumed that $S=\{0,1,\ldots,n-1\}$, for a subset $X\subseteq S$, $X=\{x_0,x_1,\ldots,x_r\}$, $x_0<x_2<\cdots<x_r$. Denoted by $\Theta_i=(x_i,x_{i+1})$, $i=0,1,\ldots,r-1$, and $\Theta_r=(x_r,n)$, take $X$ as a subset of the complete residue system $\mathcal{O} = \{0, 1, \ldots, n-1\}$, $\mod(n)$, then $\mathcal{O}$ may be decomposed into that

$$\mathcal{O} = [0, x_0] \cup \Theta_0 \cup \{x_1\} \cup \Theta_1 \cup \{x_2\} \cup \Theta_2 \cup \cdots \cup \{x_r\} \cup \Theta_r,$$

with a possible rotation, it may be assumed that $x_0 = 0$. Let $\Theta(X) = \bigcup_{0 \leq i \leq r} \Theta_i(X)$, clearly, $\Theta(X) = S \setminus X$. Denoted by $b_i = |\Theta_i|$, $0 \leq i \leq r$, and $h_X = \max_{0 \leq i \leq r} \{b_i\}$, that is, $h_X$ is the size of the largest gap of set $X$ in the circle $\mathcal{O}$.

For a non-negative integer $j$, $0 \leq j < n$, and for a subset $X \subseteq S$, define

$$\varphi_j(X) = j + \sum_{x \in X} x \quad \mod(n).$$

we set

$$\mathfrak{L}_j = \left\{ X \mid X \in \binom{S}{k}, \varphi_j(X) \in \{0, 1, \ldots, h_X - 1\} \right\}.$$

It is not difficult to demonstrate that $\mathfrak{L}_j$ is a $(k+1, k)$-Turán family: Suppose that $Z$ is a $(k+1)$-subset of $S$. If $\varphi_j(Z) = \alpha \in Z$, take $X = Z \setminus \alpha$, then $\varphi_j(X) = 0$, that is, $X \in \mathfrak{L}_j$; if $\varphi_j(Z) = \alpha \in \Theta_i(Z)$, take $X = Z \setminus \{x_i\}$, then $\varphi_j(X) = \alpha - x_i \leq b_i < h_X$, i.e. $X \in \mathfrak{L}_j$.

Besides, we know that each $k$-subset $X$ of $S$ occurs exactly $h_X$ $\mathfrak{L}_j$'s, hence

$$T(k+1, k, n) \leq \frac{1}{n} \sum_{j=0}^{n-1} |\mathfrak{L}_j| \leq \frac{1}{n} \sum_{X \in \binom{S}{k}} h_X$$

$$\leq \frac{1}{n} \times \sum_h h \cdot \lambda(s, k, h) \leq \frac{1 + \tau(k)}{k-1} \cdot C_n^k.$$

□

We have made three lists for the precise values of numbers $c_0(k), c_1(k)$ and $|\mathfrak{L}_j|$ with aid of computer, which are attached in the end of the paper as an appendix.

# Appendix

| k | $c_0(k)$ | k | $c_0(k)$ | k | $c_0(k)$ | k | $c_0(k)$ |
|---|---|---|---|---|---|---|---|
| 1 | 0.000000 | 11 | 0.186127 | 21 | 0.202151 | 31 | 0.207631 |
| 2 | 0.000000 | 12 | 0.189120 | 22 | 0.202917 | 32 | 0.208013 |
| 3 | 0.027084 | 13 | 0.191533 | 23 | 0.203605 | 33 | 0.208364 |
| 4 | 0.091329 | 14 | 0.193566 | 24 | 0.204267 | 34 | 0.208667 |
| 5 | 0.127047 | 15 | 0.195307 | 25 | 0.204856 | 35 | 0.208974 |
| 6 | 0.148383 | 16 | 0.196800 | 26 | 0.205400 | 36 | 0.209271 |
| 7 | 0.161989 | 17 | 0.198148 | 27 | 0.205932 | 37 | 0.209543 |
| 8 | 0.171132 | 18 | 0.199298 | 28 | 0.206419 | 38 | 0.209775 |
| 9 | 0.177626 | 19 | 0.200355 | 29 | 0.206828 | 39 | 0.210028 |
| 10 | 0.182418 | 20 | 0.201297 | 30 | 0.207241 | 40 | 0.210253 |

List 1.

| k | $c_1(k)$ | k | $c_1(k)$ | k | $c_1(k)$ | k | $c_1(k)$ |
|---|---|---|---|---|---|---|---|
| 1 | 0.000000 | 11 | 0.808947 | 21 | 0.802982 | 31 | 0.800903 |
| 2 | 0.875000 | 12 | 0.807892 | 22 | 0.802688 | 32 | 0.800767 |
| 3 | 0.845679 | 13 | 0.807002 | 23 | 0.802420 | 33 | 0.800640 |
| 4 | 0.832357 | 14 | 0.806243 | 24 | 0.802174 | 34 | 0.800520 |
| 5 | 0.824731 | 15 | 0.805587 | 25 | 0.801948 | 35 | 0.800407 |
| 6 | 0.819788 | 16 | 0.805014 | 26 | 0.801740 | 36 | 0.800300 |
| 7 | 0.816324 | 17 | 0.804511 | 27 | 0.801548 | 37 | 0.800199 |
| 8 | 0.813760 | 18 | 0.804064 | 28 | 0.801369 | 38 | 0.800104 |
| 9 | 0.811787 | 19 | 0.803664 | 29 | 0.801203 | 39 | 0.800013 |
| 10 | 0.810221 | 20 | 0.803306 | 30 | 0.801048 | 40 | 0.799927 |

List 2.

| k | $L_{\min}$ | $\bar{L}$ | $\bar{L}/C_n^k$ | $(1+\tau(k))/(k-1)$ |
|---|---|---|---|---|
| 2 | 345 | 352 | 1/1.409091 | 1 |
| 3 | 2802 | 2810 | 1/1.765125 | 1/1.043184 |
| 4 | 16939 | 16977 | 1/2.118160 | 1/1.410224 |
| 5 | 81414 | 81468 | 1/2.471842 | 1/1.733762 |
| 6 | 320071 | 320180 | 1/2.830258 | 1/2.029909 |
| 7 | 1052839 | 1052987 | 1/3.196484 | 1/2.307475 |
| 8 | 2943440 | 2943686 | 1/3.573173 | 1/2.571572 |
| 9 | 7077594 | 7078001 | 1/3.962814 | 1/2.825472 |
| 10 | 14769213 | 14769646 | 1/4.367893 | 1/3.071266 |

$$L_{\min} = \min_{0 \le j < n}\{|\mathfrak{L}_j|\},\ \bar{L} = \frac{1}{n}\sum_{0 \le j < n}|\mathfrak{L}_j|,\ n = 32,\ \tau(k) = \log(k) - 1 + c_1(k) - c_0(k).$$

List 3.